\documentclass[reqno,psamsfonts,11pt]{gen-j-l}
\usepackage{amssymb,amsmath,amsbsy,amsthm,esint,setspace}
\usepackage{hyperref}
\usepackage{amsrefs}
\usepackage{enumitem}

\newtheorem{theorem}{Theorem}[section]
\newtheorem{lemma}[theorem]{Lemma}
\newtheorem{proposition}[theorem]{Proposition}

\newtheorem{definition}{Definition}[section]

\newtheorem*{remark}{Remark}

\begin{document}

\title[Structure of correlations for hard spheres]
{Structure of correlations for the Boltzmann-Grad limit of hard spheres}

\author{Ryan Denlinger}

\begin{abstract}
 We consider a gas of $N$ identical hard spheres in the whole space,
and we enforce the Boltzmann-Grad scaling. We may suppose that the particles
are essentially independent of each other at some initial time; even
so, correlations will be created by the dynamics. We will prove a
structure theorem for the correlations which develop at positive time.
Our result generalizes a previous result which states that 
there are phase points where the three-particle marginal density factorizes
into two-particle and one-particle parts, while further factorization
is impossible. The result depends on uniform bounds which are known to hold
on a small time interval, or globally in time when the mean free path
is large.
\end{abstract}


\maketitle

\section{Introduction}
\label{intro}

We are interested in the problem of deriving macroscopic evolutionary
equations from a Newtonian gas of $N$ identical hard spheres, each having
diameter $\varepsilon>0$ and set in the spatial domain $\mathbb{R}^d$
for some $d\geq 2$. The formal scaling we will concern ourselves with is the
Boltzmann-Grad scaling, which means that the mean free path for a particle of
gas is of order one. Assuming that particles are initially independent of
one another, the expected evolution equation in this scaling is
Boltzmann's equation with hard sphere collision kernel.
Our goal is to refine known results on the propagation of chaos; we will study
the structure of correlations on parts of the phase space where the
pure factorization structure is necessarily destroyed.
Note carefully that we are interested in the correlations between different
particles' configurations at a \emph{fixed} time.
We have proven, in our previous work,
that on some parts
of the reduced 
phase space, the marginal density for three particles factorizes
into two-particle and one-particle contributions, while further factorization
is impossible. Our aim is to generalize that result to correlations
of $m-1$ particles, for any finite $m$.

Dynamically-induced correlations in Newtonian hard sphere gases have received
some attention in the recent literature; we will remark on two results in
particular. Pulvirenti and Simonella analyzed the size of higher-order
correlations in the context of Lanford's theorem. \cite{PS2014,L1975} 
Remarkably, the authors were able to quantify the correlations even among
$\varepsilon^{-\alpha}$ particles for some $\alpha > 0$.
This work required a sophisticated analysis of many-recollision events,
and employed a special representation formula to make clear the obstructions
to factorization. The other work is a derivation of linear hydrodynamics
by Bodineau, Gallagher and Saint-Raymond. \cite{BGSR2015III} This work relied
on a perturbative expansion accounting for corrections to
(linearized) factorization.
Contrary to one's naive expectation,
 the authors were able to
quantitatively control corrections of all orders globally in time in a
weighted $L^2$ norm.

Morally speaking,
the results we prove in this work  ought to
 show that, conditional on certain
$L^\infty$ estimates and the factorization of the initial data, the
$s$th marginal reduces to a tensor product of $(s-1)$st marginal and the
first marginal, as long as the backwards trajectory of
\emph{one} particle is free. The reason we cannot remove the words
``morally'' and ``ought to'' is that the result requires the deletion
of certain explicit
sets upon which
the convergence either may not or does not hold.  In any case, if we view the first
$(s-1)$ particle configurations as fixed and choose the $s$th particle's
configuration ``randomly,'' it is highly unlikely we will land on the
exceptional set; this is true even if the first $(s-1)$ particles
have very complicated backwards trajectories. For this reason, the
deletion of an exceptional set does not alter significantly the
interpretation of our results.

We have drawn particular inspiration
from the methods of Pulvirenti and Simonella \cite{PS2014}, 
and we have explicitly used their idea
of employing an intermediate Boltzmann-Enskog-type hierarchy between the
BBGKY and Boltzmann hierarchies. However, we have avoided using their 
special
representation formula for correlations. Instead, we employ an
unsymmetric Boltzmann-Enskog hierarchy (defined in our previous work
\cite{De2017}) which tracks correlations between
just the first $m-1$
 particles. This is convenient because we can show that
the intermediate
 hierarchy propagates partial factorization in an \emph{exact} sense.
The theorem then follows by quantitatively comparing the
BBGKY and Boltzmann-Enskog pseudo-dynamics; our primary concern is
 making such stability estimates precise.

We rely heavily on the developments of our previous work \cite{De2017}; indeed,
all we do is refine part \emph{(ii)} of 
Theorem 2.1 from that work to include the case of arbitrarily many
correlated particles (instead of just two correlated particles). For
this reason, we will only briefly summarize the previous results
 that are relevant here, and fill in the missing ideas and
estimates that are needed to prove our main theorem.

\subsection*{Organization.} The notation and main theorem
are given in Section \ref{sec:2}. The BBGKY hierarchy, with basic
results, is recalled in Section \ref{sec:4}. In Section \ref{sec:5},
we recall an unsymmetric Boltzmann-Enskog
hierarchy from Appendix A of our previous work. \cite{De2017}
A crucial stability result is proven in
Section \ref{sec:6}; the remainder of the convergence proof proceeds
as in the aforementioned work. \cite{De2017}

\section{Notation and Main Results}
\label{sec:2}

We consider $N$ identical hard spheres with diameter $\varepsilon > 0$,
 positions $x_1,x_2,\dots,x_N \in
\mathbb{R}^d$, and velocities
$v_1,v_2,\dots,v_N \in \mathbb{R}^d$. The tuple of all particle
positions is written $X_N = (x_1,x_2,\dots,x_N) \in \mathbb{R}^{dN}$,
and the tuple of all particle velocities is written
$V_N = (v_1,v_2,\dots,v_N) \in \mathbb{R}^{dN}$. We also write
$Z_N = (z_1,z_2,\dots,z_N) = (X_N,V_N) \in \mathbb{R}^{2dN}$.
The Boltzmann-Grad scaling $N\varepsilon^{d-1}=\ell^{-1}$, for fixed
$\ell > 0$, is assumed throughout.
The $N$-particle phase space is the set
\begin{equation}
\label{eq:s2-D-N}
\mathcal{D}_N =  \left\{\left. Z_N = (X_N,V_N) \in \mathbb{R}^{2dN} \right|
\forall 1\leq i < j \leq N, \; |x_i-x_j|>\varepsilon \right\} 
\end{equation}
As long as $Z_N \in \mathcal{D}_N$ we allow particles to move in
straight lines with constant velocity: $\dot{X}_N = V_N$, 
$\dot{V}_N = 0$. Specular reflection is enforced at the
boundary $\partial \mathcal{D}_N$; up to deletion of a zero measure
set, all collisions are binary, non-grazing, and linearly ordered
in time. If the $i$th and $j$th particles collide at time $t_0$ with
$x_j (t_0^-)  = x_i (t_0^-) + \varepsilon \omega$ then the velocities
transform according to the following rule:
\begin{equation}
\begin{aligned}
v_i (t_0^+) & = v_i (t_0^-) + \omega \omega \cdot \left(
v_j (t_0^-) - v_i (t_0^-) \right) \\
v_j (t_0^+) & = v_j (t_0^-) - \omega\omega \cdot \left(
v_j (t_0^-) - v_i (t_0^-) \right)
\end{aligned}
\end{equation}
The collective flow of $N$ identical hard spheres of diameter
$\varepsilon > 0$ defines a measurable map 
$\psi_N^t : \overline{\mathcal{D}_N} \rightarrow \overline{\mathcal{D}_N}$
preserving the
Lebesgue measure on $\mathcal{D}_N$.
We define a measurable involution $Z_N \mapsto Z_N^*$ on
$\partial \mathcal{D}_N$ which is defined almost everywhere by
the following properties:
\begin{equation*}
(\textnormal{a.e. } Z_N \in \partial \mathcal{D}_N) \qquad
\left( \lim_{t\rightarrow 0^-} \psi_N^{t} Z_N \right)^*
= \lim_{t \rightarrow 0^+} \psi_N^t Z_N
\end{equation*}
\begin{equation*}
(\textnormal{a.e. } Z_N \in \partial \mathcal{D}_N) \qquad
\left( \lim_{t\rightarrow 0^+} \psi_N^{t} Z_N \right)^*
= \lim_{t \rightarrow 0^-} \psi_N^t Z_N
\end{equation*}

We introduce a probability measure $f_N (0)$ on $\mathcal{D}_N$ and
define $f_N (t)$ to be the pushforward of $f_N (0)$ under 
$\psi_N^t$. Since $\psi_N^t$ preserves the Lebesgue measure on
$\mathcal{D}_N$, this says that
\begin{equation}
f_N (t,Z_N) = f_N \left(0,\psi_N^{-t} Z_N\right)
\end{equation}
We assume that $f_N (0)$ is symmetric under interchange of particle
indices; since the particles are identical, it follows that
$f_N (t)$ is symmetric as well. We extend $f_N (t)$ by zero to be
defined on all of $\mathbb{R}^{2dN}$.

We define the marginals $f_N^{(s)} (t)$, $1\leq s \leq N$, by the
formula
\begin{equation}
f_N^{(s)} (t,Z_s) = \int_{\mathbb{R}^{2d(N-s)}}
f_N (t,Z_N) dz_{s+1} \dots dz_N
\end{equation}
Then the support of $f_N^{(s)} (t)$ is contained in the closure of
$\mathcal{D}_s$, where
\begin{equation}
\mathcal{D}_s = \left\{ \left. Z_s = (X_s,V_s) \in \mathbb{R}^{2ds}
\right| \forall 1\leq i < j \leq s,\; |x_i - x_j|>\varepsilon
\right\}
\end{equation}
Note carefully that $\mathcal{D}_s$ depends on $\varepsilon$ whenever
$s\geq 2$; however, this dependence is suppressed in our notation.
The flow of $s$ identical hard spheres of diameter $\varepsilon$ is
written $\psi_s^t : \mathcal{D}_s \rightarrow \mathcal{D}_s$; again, the
implicit dependence on $\varepsilon$ is suppressed in the notation.
We also define
$E_s (Z_s) = \frac{1}{2} \sum_{i=1}^s |v_i|^2$, which is the total
energy of $s$ particles.

In order to state our main results, we will require a
 notion of \emph{deletion} of particles; this
will be helpful in defining the exceptional set where convergence
may fail. For any $1\leq k \leq s$ we define
\begin{equation}
r_{k} Z_s = \left( z_1,\dots,z_{k-1},z_{k+1},\dots,z_s\right)
\end{equation}
In other words, if $Z_s$ is any \emph{ordered list} of particle
configurations, then $r_k Z_s$ is the same list with the $k$th entry
removed. For example,
\begin{equation}
r_4 r_3 Z_5 = r_4 \left(z_1,z_2,z_4,z_5\right) = \left(z_1,z_2,z_4\right)
\end{equation}
In this example, $r_4$ deletes the fourth particle in the list,
\emph{not} the particle with initial label $z_4$.

\begin{remark}
An alternative notation for particle deletion
would be possible if we chose to associate
with each particle a label, so that $Z_2$ really denotes
$\left\{ (z_1,\textnormal{``1''}),
(z_2,\textnormal{``2''})\right\}$. Then we could define $r_2$ to be
the deletion of the particle with \emph{label} equal to $2$.
However, such notation is
not needed here, so instead $Z_s$ is to be viewed simply as
an ordered list of points in $\mathbb{R}^{2d}$.
\end{remark}

Now for any $Z_s \in \mathcal{D}_s$ we define the \emph{set} of
points $\mathcal{J}_s Z_s \subset \mathbb{R}^{2d}$ as follows:
\begin{equation}
\mathcal{J}_s Z_s = 
\bigcup_{
\substack{
\tau_2,\dots,\tau_s \geq 0 \\
1\leq k_j \leq j
}}
\psi_1^{\tau_2 + \dots + \tau_s}
r_{k_2} \psi_2^{-\tau_2} r_{k_3} \psi_3^{-\tau_3} \dots
r_{k_{s}} \psi_s^{-\tau_s} Z_s
\end{equation}
Note that $\mathcal{J}_s Z_s$ is a finite set, and its cardinality
can even be controlled in terms of $s$. \cite{Va1979,I1989,BFK1998}
We want to view the first $m-1$ particles as ``interacting'' and
the remaining $s-m+1$ particles as ``free.'' Hence we define:
\begin{equation}
\label{eq:s2-G-s-m}
\mathcal{G}_{s|m} = \left\{ Z_s  \in 
\overline{\mathcal{D}_s} \left|
\begin{aligned}
& \forall m \leq i \leq s,\;
\forall \tau > 0, \; \forall (x^0,v^0) \in \mathcal{J}_{m-1} Z_{m-1}, \\
& \qquad \qquad\qquad 
 |(x_i - x^0) - (v_i - v^0)\tau| \geq \varepsilon \\
& \textnormal{and} \; \forall m\leq i \neq j \leq s,\; \forall
\tau > 0, \\
& \qquad \qquad \qquad
|(x_i - x_j) - (v_i-v_j)\tau| \geq \varepsilon 
\end{aligned}
\right. \right\}
\end{equation}
The condition $Z_s \in \mathcal{G}_{s|m}$ means that the last
$s-m+1$ particles are free under the backwards flow
\emph{no matter} the history of the first $m-1$ particles, 
\emph{including} the possibility that some of the first $m-1$
particles may be ``removed'' from the interaction at arbitrary
intermediate times. We warn the
reader that this is only a heuristic explanation and the true
definition is given by (\ref{eq:s2-G-s-m}).

We will also need a condition which forces
particles to disperse from one another. Hence for any
$\eta > 0$ we define:
\begin{equation}
\hat{\mathcal{U}}_s^\eta = \left\{ Z_s \in \overline{\mathcal{D}}_s
\left| 
\begin{aligned}
& \forall (x^0,v^0),(x^1,v^1)\in \mathcal{J}_s Z_s \; : \;
(x^0,v^0) \neq (x^1,v^1), \\
& \qquad \qquad \qquad \qquad \qquad
|v^0 - v^1|>\eta
\end{aligned}
\right. \right\}
\end{equation}
We also define
\begin{equation}
\label{eq:s2-K-s}
\mathcal{K}_s = \left\{ Z_s \in \overline{\mathcal{D}_s} \left|
\forall \tau > 0,\; \psi_s^{-\tau} Z_s = (X_s - V_s \tau,V_s)
\right. \right\}
\end{equation}
\begin{equation}
\label{eq:s2-U-s-eta}
\mathcal{U}_s^\eta = \left\{ Z_s \in \overline{\mathcal{D}_s}
 \left| \inf_{1\leq i < j \leq s} |v_i - v_j| > \eta \right. \right\}
\end{equation}

\begin{definition}
\label{def:s2-chaos}
Let us be given, for each $N\in\mathbb{N}$, a sequence of densities
$F_N = \left\{ f_N^{(s)} \right\}_{1\leq s \leq N}$,
with each $f_N^{(s)}$ defined on $\mathcal{D}_s$ and symmetric with
respect to particle interchange. Then
$\left\{ F_N \right\}_{N\in\mathbb{N}}$ is 
\emph{$(m-1)$-nonuniformly $f$-chaotic} for some density $f(z)$  if, for
some $\kappa \in (0,1)$, there holds for each integer
$3 \leq m^\prime \leq m$, every $s \geq m^\prime-1$, and
all $R>0$ that
\begin{equation}
\begin{aligned}
& \lim_{N\rightarrow \infty}
\left\Vert \left( f_N^{(s)}  - f_N^{(m^\prime-1)}\otimes 
f^{\otimes (s-m^\prime+1)} \right) (Z_s)
\mathbf{1}_{Z_s \in \mathcal{G}_{s|m^\prime}
 \cap \hat{\mathcal{U}}_s^{\eta(\varepsilon)}}
\mathbf{1}_{\frac{1}{2} \sum_{i=1}^s |v_i|^2 \leq R^2} 
\right\Vert_{L^\infty_{Z_s}} \\
& \qquad \qquad \qquad \qquad \qquad \qquad \qquad \qquad
\qquad \qquad \qquad \qquad \qquad \qquad   = 0
\end{aligned}
\end{equation}
and, for each integer $s\geq 1$ and all $R>0$,
\begin{equation}
\begin{aligned}
& \lim_{N\rightarrow \infty}
\left\Vert \left( f_N^{(s)} (Z_s) -  f^{\otimes s} 
(Z_s)\right)
\mathbf{1}_{Z_s \in \mathcal{K}_s \cap \mathcal{U}_s^{\eta(\varepsilon)}}
\mathbf{1}_{\frac{1}{2} \sum_{i=1}^s |v_i|^2 \leq R^2} 
\right\Vert_{L^\infty_{Z_s}} = 0
\end{aligned}
\end{equation}
where $\eta (\varepsilon) = \varepsilon^\kappa$.
If $\left\{ F_N \right\}_{N\in \mathbb{N}}$ is $(m-1)$-nonuniformly
$f$-chaotic for every $3\leq m \in \mathbb{N}$ then we say that
$\left\{ F_N \right\}_{N\in\mathbb{N}}$ is
$\infty$-nonuniformly $f$-chaotic.
\end{definition}

\begin{remark}
Note that $L^\infty_{Z_s}$ refers to the essential supremum norm,
since the marginals are only defined up to sets of measure zero.
\end{remark}

\begin{remark}
The definition of $(m-1)$-nonuniform chaoticity is not exactly the
same as the notion of $2$-nonuniform chaoticity we have introduced
previously \cite{De2017} when $m=3$. Nevertheless, the two notions are almost
the same, in terms of the complexity of sets involved in the definition.
\end{remark}

Recall the Boltzmann equation for hard spheres,
\begin{equation}
\label{eq:s2-boltz}
\left( \partial_t + v\cdot \nabla_x\right) f(t) =
\ell^{-1} Q(f(t),f(t)) 
\end{equation}
where
\begin{equation}
\label{eq:s2-coll1}
Q(f,f) = \int_{\mathbb{R}^d \times \mathbb{S}^{d-1}}
\left[ \omega \cdot (v_1 - v)\right]_+ \left(
f(x,v^*) f(x,v_1^*) - f(x,v) f(x,v_1)\right) d\omega dv_1
\end{equation}
We are now ready to state our main result.
\begin{theorem}
\label{thm:s2-chaos}
For each $N\in\mathbb{N}$, let $F_N (t) = \left\{ f_N^{(s)} (t)
\right\}_{1\leq s\leq N}$ solve the hard sphere BBGKY
hierarchy, enforcing the Boltzmann-Grad scaling $N\varepsilon^{d-1}
=\ell^{-1}$. Assume that each $f_N^{(s)} (t)$ is symmetric with respect
to particle interchange. Let $f (t,x,v)$ solve the Boltzmann equation 
(\ref{eq:s2-boltz}) for $0\leq t \leq T$; furthermore, assume that
$f(t) \in W^{1,\infty} (\mathbb{R}^d \times \mathbb{R}^d)$, $f(t) \geq 0$, 
$\int f(t) dx dv = 1$, and that there exists $\beta_T > 0$ such that
\begin{equation}
\sup_{0\leq t \leq T} \sup_{x,v} e^{\frac{1}{2} \beta_T |v|^2}
f(t,x,v) < \infty
\end{equation}
Further suppose that there exists $\tilde{\beta}_T > 0$, 
$\tilde{\mu}_T \in \mathbb{R}$ such that
\begin{equation}
\sup_{N\in\mathbb{N}} 
\sup_{1\leq s \leq N} \sup_{0\leq t\leq T} \sup_{Z_s \in \mathcal{D}_s}
e^{\tilde{\beta}_T E_s (Z_s)} e^{\tilde{\mu}_T s}
\left| f_N^{(s)} (t,Z_s) \right| < \infty
\end{equation}
Then we have the following: \\
(i) If $\left\{ F_N (0) \right\}_{N\in\mathbb{N}}$ is
$(m-1)$-nonuniformly $f(0)$-chaotic, then for each $t\in [0,T]$,
$\left\{ F_N (t) \right\}_{N\in\mathbb{N}}$ is
$(m-1)$-nonuniformly $f(t)$-chaotic (with the same $\kappa$). \\
(ii) If
$\left\{ F_N (0) \right\}_{N\in\mathbb{N}}$ is
$\infty$-nonuniformly $f(0)$-chaotic, then for each $t\in [0,T]$,
$\left\{ F_N (t) \right\}_{N\in \mathbb{N}}$ is
$\infty$-nonuniformly $f(t)$-chaotic (with the same $\kappa$).
\end{theorem}

\begin{remark}
We have stated Theorem \ref{thm:s2-chaos} without any explicit error
estimates for simplicity. However, it is not hard to extract
quantitative estimates from the proof. Note that good  error
estimates cannot
be expected even for $s \approx \log N$ due to our reliance on wildly
divergent bounds on the number of collisions
of hard spheres. \cite{BFK1998}
\end{remark}

\section{The BBGKY Hierarchy}
\label{sec:4}

The marginals $f_N^{(s)} (t)$ solve a set of equations called the
BBGKY hierarchy (Bogoliubov-Born-Green-Kirkwood-Yvon). The hierarchy
is written as follows, for $1\leq s < N$ (the $s=N$ component just
obeys Liouville's equation):
\begin{equation}
\left( \partial_t + V_s \cdot \nabla_{X_s} \right) f_N^{(s)} (t,Z_s)
= (N-s) \varepsilon^{d-1} C_{s+1} f_N^{(s+1)} (t,Z_s)
\end{equation} 
The specular boundary condition
$f_N^{(s)} (t,Z_s^*) = f_N^{(s)} (t,Z_s)$ is enforced along
$\partial \mathcal{D}_s$; it actually holds for
$\textnormal{a.e. } (t,Z_s) \in [0,T]\times
\partial \mathcal{D}_s$ ($T>0$ is arbitrary).
The \emph{collision operator} $C_{s+1}$ is written
\begin{equation}
C_{s+1} f_N^{(s+1)} (Z_s) = C_{s+1}^+ f_N^{(s+1)} (Z_s) -
C_{s+1}^- f_N^{(s+1)} (Z_s)
\end{equation}
\begin{equation}
C_{s+1}^{\pm} = \sum_{i=1}^s C_{i,s+1}^{\pm}
\end{equation}
\begin{equation}
\begin{aligned}
& C_{i,s+1}^- f_N^{(s+1)} (t,Z_s)
 = \int_{\mathbb{R}^d \times \mathbb{S}^{d-1}}
d\omega dv_{s+1}
\left[ \omega \cdot (v_{s+1}-v_i)\right]_-  \times \\
& \qquad \qquad \qquad\qquad \qquad \qquad \qquad \qquad \times
f_N^{(s+1)} \left( t,Z_s,x_i+\varepsilon \omega,v_{s+1}\right)
\end{aligned}
\end{equation}
\begin{equation}
\begin{aligned}
& C_{i,s+1}^+ f_N^{(s+1)} (t,Z_s)
= \int_{\mathbb{R}^d \times \mathbb{S}^{d-1}}
d\omega dv_{s+1} 
\left[ \omega \cdot (v_{s+1}-v_i)\right]_+  \times \\
& \qquad \qquad \qquad \qquad \qquad \qquad
\times f_N^{(s+1)} \left( t,\dots,x_i,v_i^*,\dots,
x_i+\varepsilon \omega,v_{s+1}^*\right)
\end{aligned}
\end{equation}
Here $v_i^* = v_i + \omega \omega \cdot (v_{s+1}-v_i)$ and
$v_{s+1}^* = v_{s+1}-\omega \omega \cdot (v_{s+1}-v_i)$.

The BBGKY hierarchy is well-posed locally in time in suitable
$L^\infty$ norms. \cite{L1975,GSRT2014}
 Roughly speaking, as long as the
initial data $\left\{ f_N^{(s)} (0) \right\}_{1\leq s \leq N}$ satisfies
the following bound, for some $\beta_0 > 0$, $\mu_0 \in \mathbb{R}$,
\begin{equation}
\sup_{1\leq s \leq N} \sup_{Z_s \in \mathcal{D}_s}
\left| f_N^{(s)} (0,Z_s)\right| e^{\beta_0 E_s (Z_s)}
e^{\mu_0 s} \leq 1
\end{equation}
then, in the Boltzmann-Grad scaling $N\varepsilon^{d-1} = \ell^{-1}$,
on a small time interval $T_L < C_d \ell e^{\mu_0} 
\beta_0^{\frac{d+1}{2}}$, the BBGKY hierarchy has a unique solution
satisfying the bound
\begin{equation}
\label{eq:s4-unif-est}
\sup_{0 \leq t \leq T_L}
\sup_{1\leq s \leq N} \sup_{Z_s \in \mathcal{D}_s}
\left| f_N^{(s)} (t,Z_s)\right| e^{\frac{1}{2} \beta_0 E_s (Z_s)}
e^{(\mu_0-1) s} \leq 1
\end{equation}
The well-posedness statement can be made more precise by using
time-dependent weights but this is not necessary for any of our
results. In fact the estimate (\ref{eq:s4-unif-est}) is sufficient
to guarantee that the formal series we write are bounded, uniformly
in $N$, for a short time. Our arguments can be iterated in time
 for as long as uniform bounds are available;
this is how we ultimately deduce Theorem \ref{thm:s2-chaos}. 

Let us define the operators $T_s (t)$ which act on functions
$f^{(s)}  : \mathcal{D}_s \rightarrow \mathbb{R}$ as follows:
\begin{equation}
\left( T_s (t) f^{(s)} \right) (Z_s) = f^{(s)} \left(\psi_s^{-t} Z_s\right)
\end{equation}
Then the functions $f_N^{(s)} (t)$ solve the following \emph{mild form}
of the BBGKY hierarchy:
\begin{equation}
f_N^{(s)} (t) = T_s (t) f_N^{(s)} (0) + (N-s)\varepsilon^{d-1} 
\int_0^t T_s (t-t_1) C_{s+1} f_N^{(s+1)} (t_1) dt_1
\end{equation}
Iterating this formula in the standard way \cite{L1975,GSRT2014},
 we express $f_N^{(s)} (t)$ as a finite sum
of terms involving only the initial data:
\begin{equation}
\begin{aligned}
\label{eq:s4-duhamel-ser}
& f_N^{(s)} (t) = \sum_{k=0}^{N-s} a_{N,k,s} \times \\
& \qquad  \times \int_0^t \int_0^{t_1} \dots
\int_0^{t_{k-1}} T_s (t-t_1) C_{s+1} \dots
T_{s+k} (t_k) f_N^{(s+k)} (0) dt_k \dots dt_1
\end{aligned}
\end{equation}
where
\begin{equation}
a_{N,k,s} = \frac{(N-s)!}{(N-s-k)!} \varepsilon^{k(d-1)}
\end{equation}

Following the arguments of Lanford \cite{L1975,GSRT2014,De2017}, we can define \emph{pseudo-trajectories}
which encode the possible (re)collision sequences which contribute to
the solution $f_N^{(s)} (t)$ of the BBGKY hierarchy. Given a final
state $Z_s \in \mathcal{D}_s$, along with times $0 \leq t_k \leq \dots \leq
t_1 \leq t$, velocities $v_{s+1},\dots,v_{s+k} \in \mathbb{R}^{d}$,
impact parameters $\omega_1,\dots,\omega_k \in \mathbb{S}^{d-1}$, and
indices $i_j \in \left\{ 1,\dots,s+j-1\right\}$, we define the point
\begin{equation}
Z_{s,s+k} \left[ Z_s,t;t_1,\dots,t_k;v_{s+1},\dots,v_{s+k};
\omega_1,\dots,\omega_k;i_1,\dots,i_k \right] \in 
\mathcal{D}_{s+k}
\end{equation}
We think of $Z_{s,s+k}$ as being the image of $Z_s$ under a sequence
of $k$ particle creations at times $t_j$. We evolve backwards the
point $Z_s$ under the hard sphere flow for a time $t-t_1$; then,
we create a particle adjacent to particle $i_1$,
so $x_{s+1} = x_{i_1} + \varepsilon \omega_1$; we force a collisional
change of variables, if needed, to place all particles in
a \emph{pre-collisional} state. We then continue the backwards flow
of $s+1$ particles for a time $t_1-t_2$, then create another particle, and
so forth.

We define iterated collision kernels
\begin{equation}
b_{s,s+k} \left[ Z_s,t; \left\{ t_j,v_{s+j},\omega_j,i_j
\right\}_{j=1}^k \right]
\end{equation}
which simply records the accumulated product of impact
parameters, e.g. $\omega \cdot (v_1 - v_2)$
for a collision involving the particles $1$ and $2$;
 finally, the iterated Duhamel series
(\ref{eq:s4-duhamel-ser}) becomes
\begin{equation}
\label{eq:s4-duhamel-dyn}
\begin{aligned}
& f_N^{(s)} (t,Z_s) = \sum_{k=0}^{N-s} a_{N,k,s} \times \\
& \times \sum_{i_1 = 1}^s \dots \sum_{i_k = 1}^{s+k-1}
\int_0^t \dots \int_0^{t_{k-1}} \int_{\mathbb{R}^{dk}}
\int_{(\mathbb{S}^{d-1})^k} \left(
\prod_{m=1}^k d\omega_m dv_{s+m} dt_m \right) \times \\
& \times \left( b_{s,s+k} \left[ \cdot \right] 
f_N^{(s+k)} \left(0,Z_{s,s+k}\left[\cdot\right]\right)
\right) \left[ Z_s,t; \left\{ t_j,v_{s+j},\omega_j,i_j
\right\}_{j=1}^k \right]
\end{aligned}
\end{equation}
We refer to Section 7 of our previous work \cite{De2017} for further details.

\section{The Unsymmetric Boltzmann-Enskog Hierarchy}
\label{sec:5}

We summarize the developments of Appendix A, sections A.1 and
A.2, of our previous work. \cite{De2017} We will not quote precise results; the
reader may refer to Appendix A of that work for theorems and
proofs.

We are going to write down an unsymmetric Boltzmann-Enskog-type
hierarchy which tracks correlations between the first $m-1$
particles. Let us define the unsymmetric $s$-particle phase space,
where $m\geq 2$ is fixed and $s\geq m-1$:
\begin{equation}
\tilde{\mathcal{D}}_s = \left\{ \left.
Z_s = (X_s,V_s) \in \mathbb{R}^{2ds} \right| \;
\forall 1\leq i < j \leq m-1,\; |x_i-x_j|>\varepsilon \right\}
\end{equation} 
Furthermore, define the collision operators
\begin{equation}
\tilde{C}_{s+1} = \sum_{i=1}^s \left(
\tilde{C}_{i,s+1}^+ - \tilde{C}_{i,s+1}^- \right)
\end{equation}
where
\begin{equation}
\begin{aligned}
& \tilde{C}_{i,s+1}^- g_\varepsilon^{(s+1)} (t,Z_s) =
\int_{\mathbb{R}^d \times \mathbb{S}^{d-1}} d\omega dv_{s+1} \left[
\omega \cdot (v_{s+1}-v_i )\right]_- \times \\
& \qquad \qquad \qquad \qquad \qquad \times g_\varepsilon^{(s+1)}
(t,\dots,x_i,v_i,\dots,x_i+\varepsilon \omega,v_{s+1})
\end{aligned}
\end{equation}
\begin{equation}
\begin{aligned}
& \tilde{C}_{i,s+1}^+ g_\varepsilon^{(s+1)} (t,Z_s) =
\int_{\mathbb{R}^d \times \mathbb{S}^{d-1}} d\omega dv_{s+1} \left[
\omega \cdot (v_{s+1}-v_i) \right]_+ \times \\
& \qquad \qquad\qquad \qquad\qquad \times g_\varepsilon^{(s+1)}
(t,\dots,x_i,v_i^*,\dots,x_i+\varepsilon \omega,v_{s+1}^*)
\end{aligned}
\end{equation}
with $v_i^* = v_i + \omega \omega \cdot (v_{s+1}-v_i)$ and 
$v_{s+1}^* = v_{s+1}-\omega \omega \cdot (v_{s+1}-v_i)$.
The unsymmetric Boltzmann-Enskog hierarchy is then written, for
$Z_s \in \tilde{\mathcal{D}}_s$, $s\geq m-1$,
\begin{equation}
\label{eq:s5-UBEH}
\left( \partial_t + V_s \cdot \nabla_{X_s}\right)
g_\varepsilon^{(s)} = \ell^{-1} 
\tilde{C}_{s+1} g_\varepsilon^{(s+1)} (t) \qquad
\textnormal{(if } s \geq m-1 \textnormal{)}
\end{equation}
with boundary condition
\begin{equation}
g_\varepsilon^{(s)} (t,Z_s^*) = g_\varepsilon^{(s)} (t,Z_s) \qquad
\textnormal{a.e. }  (t,Z_s) \in 
[0,T)\times \partial \tilde{\mathcal{D}}_s
\end{equation}
and initial conditions $g_\varepsilon^{(s)} (0,Z_s)$ given for
$s\geq m-1$ and $Z_s \in \overline{\mathcal{D}}_s$. We also define
the function $g_\varepsilon (t,x,v)$, $t \in [0,T],\;
x,v \in \mathbb{R}^d$, to be the solution of the equation
\begin{equation}
\label{eq:s5-BEE}
\left( \partial_t + v\cdot \nabla_x\right) g_\varepsilon (t) =
\ell^{-1} \tilde{C}_2 \left( g_\varepsilon (t) \otimes
g_\varepsilon (t) \right)
\end{equation}
with given initial data $g_\varepsilon (0)$.

The unsymmetric Boltzmann-Enskog hierarchy (\ref{eq:s5-UBEH}) and
the Boltzmann-Enskog equation (\ref{eq:s5-BEE}) are locally well-posed
when the data is in appropriate weighted $L^\infty$ spaces, just like
the BBGKY hierarchy. The proof proceeds as in the proof of
Lanford's theorem. Another important result is that the unsymmetric
Boltzmann-Enskog hierarchy propagates \emph{partial factorization},
in the following sense: Suppose that for all $s\geq m-1$ we have
\begin{equation}
g_\varepsilon^{(s)} (0) = g_\varepsilon^{(m-1)} (0) \otimes
g_\varepsilon (0)^{\otimes (s-m+1)}
\end{equation}
along with weighted $L^\infty$ bounds at the initial time. Then on
a small time interval $0\leq t \leq T$ we also have
\begin{equation}
g_\varepsilon^{(s)} (t) = g_\varepsilon^{(m-1)} (t) \otimes
g_\varepsilon (t)^{\otimes (s-m+1)}
\end{equation}
for all $s\geq m-1$. To prove this, one constructs a solution which
satisfies the partial factorization \emph{ansatz}, then the conclusion
follows by uniqueness. This is similar to the proof that the
Boltzmann hierarchy propagates factorization.

To conclude, we remark that the unsymmetric Boltzmann-Enskog hierarchy
has associated pseudo-trajectories just like the BBGKY hierarchy,
which we denote
\begin{equation}
\tilde{Z}_{s,s+k} \left[ Z_s,t;
\left\{ t_j,v_{s+j},\omega_j,i_j\right\}_{j=1}^k \right]
\end{equation}
We denote the associated iterated collision kernel by
\begin{equation}
\tilde{b}_{s,s+k} \left[ Z_s,t;\left\{ t_j,v_{s+j},\omega_j,i_j
\right\}_{j=1}^k \right]
\end{equation}

\begin{remark}
Psuedo-trajectories for the unsymmetric Boltzmann-Enskog
hierarchy are similar to pseudo-trajectories for the BBGKY or
Boltzmann hierarchies. In fact for the unsymmetric case we allow
recollisions involving the first $m-1$ particles among each
other. On the other hand, if two particles ``collide'' and at least
one of them has index $j\geq m$ then the two particles simply
pass through each other without interacting. Thus, to the extent that
the first $m-1$ particles are isolated from the remaining particles,
it follows that the
unsymmetric Boltzmann-Enskog hierarchy tracks correlations for
a cluster of $m-1$ particles.
\end{remark}

The solution $g_\varepsilon^{(s)} (t)$, $s\geq m-1$, of the
unsymmetric Boltzmann-Enskog hierarchy has the following
representation in terms of the data:
\begin{equation}
\label{eq:s5-duhamel-dyn}
\begin{aligned}
& g_\varepsilon^{(s)} (t,Z_s) = \sum_{k=0}^\infty \ell^{-k} \times \\
& \times \sum_{i_1 = 1}^s \dots \sum_{i_k = 1}^{s+k-1}
\int_0^t \dots \int_0^{t_{k-1}} \int_{\mathbb{R}^{dk}}
\int_{(\mathbb{S}^{d-1})^k} \left(
\prod_{m=1}^k d\omega_m dv_{s+m} dt_m \right) \times \\
& \times \left( \tilde{b}_{s,s+k} \left[ \cdot \right] 
g_\varepsilon^{(s+k)} \left(0,\tilde{Z}_{s,s+k}\left[\cdot\right]\right)
\right) \left[ Z_s,t; \left\{ t_j,v_{s+j},\omega_j,i_j
\right\}_{j=1}^k \right]
\end{aligned}
\end{equation}

\section{Stability of Pseudo-Trajectories}
\label{sec:6}

We finally turn to the main new estimate which allows us to conclude
Theorem \ref{thm:s2-chaos}. Both the statement and the proof largely
follow Proposition 8.3 and Proposition A.3 of our previous work \cite{De2017}.
Once we have Proposition
\ref{prop:s6-stability}, to be proven momentarily, it is a simple matter
to prove Theorem \ref{thm:s2-chaos} by estimating errors pointwise,
as in Section 12 of our previous work. \cite{De2017}
Note that Proposition \ref{prop:s6-stability} really only holds for
hard spheres because it relies on certain bounds on the number of
collisions. \cite{BFK1998}

We recall two useful lemmas from a previous work. \cite{De2017}
Combining these two lemmas, one deduces that pathological collision
events (recollisions) occur with small probability. These lemmas are
related to the collisional change of variables
$Z_s \mapsto Z_s^*$, which is interpreted as a reflection in
a suitable frame of reference.  

\begin{lemma}
\label{lemma:l1}
Fix $v \in \mathbb{R}^d \backslash \left\{ 0 \right\}$  ($d\geq 2$)
and
let $\mathbb{S}^{d-1} = \left\{ w \in \mathbb{R}^d
\left| |w| = 1 \right. \right\}$.
For any $\omega \in \mathbb{S}^{d-1}$ define
\begin{equation}
u_\omega = |v|^{-1} \left( 2 \omega \omega \cdot v - v \right)
\in \mathbb{S}^{d-1}
\end{equation}
If $\mathbb{S}_v^{d-1} =
\left\{ \omega \in \mathbb{S}^{d-1} \left|
\omega \cdot v > 0 \right. \right\}$ then the
map $\omega \mapsto u_\omega$ restricts to a diffeomorphism
$\mathbb{S}_v^{d-1} \rightarrow
\mathbb{S}^{d-1} \backslash 
\left\{ -|v|^{-1} v\right\}$.
\end{lemma}

\begin{lemma}
\label{lemma:l2}
Let $L \subset \mathbb{R}^d$ $(d\geq 2)$ be a line,
and let
$\mathbb{S}^{d-1} = \left\{
w \in \mathbb{R}^d \left| |w|=1 \right. \right\}$.
For $\rho > 0$ define the solid cylinder
\begin{equation}
\mathcal{C}_\rho = \left\{ u \in \mathbb{R}^d
\; \left| \; 
\textnormal{dist} (u,L) \leq \rho \right. \right\}
\end{equation}
Then
\begin{equation}
\int_{\mathbb{S}^{d-1}}
\mathbf{1}_{\omega \in \mathcal{C}_\rho}
d\omega \leq C_d \rho^{(d-1)/2}
\end{equation}
for some constant $C_d$ depending only on the
ambient dimension $d$ (but not depending on $\rho$ or $L$).
\end{lemma}

\begin{proposition}
\label{prop:s6-stability}
There is a constant $c_d > 0$ such that all the following holds: Assume that
\begin{equation}
\begin{aligned}
& Z_{s,s+k} \left[ Z_s,t; t_1,\dots,t_k;v_{s+1},\dots,v_{s+k};\omega_1,\dots,
\omega_k;i_1,\dots,i_k \right] = \\
& \qquad \qquad \qquad \qquad \qquad \qquad \qquad
= (X_{s+k}^\prime, V_{s+k}^\prime) \in 
\mathcal{G}_{(s+k)|m} \cap \hat{\mathcal{U}}_{s+k}^\eta
\end{aligned}
\end{equation}
and $E_{s+k} (Z_{s+k}^\prime ) \leq 2 R^2$ with $\eta < R$; then, \\
(i) for all $\tau \geq 0$ we have
\begin{equation}
\begin{aligned}
& Z_{s,s+k} \left[ Z_s,t+\tau; 
t_1+\tau,\dots,t_k+\tau;v_{s+1},\dots,v_{s+k};\omega_1,\dots,
\omega_k;i_1,\dots,i_k \right] = \\
& \qquad \qquad \qquad \qquad \qquad \qquad \qquad
= (X_{s+k}^\prime, V_{s+k}^\prime) \in 
\mathcal{G}_{(s+k)|m} \cap \hat{\mathcal{U}}_{s+k}^\eta
\end{aligned}
\end{equation}
(ii) for any $i_{k+1} \in \left\{ 1,2,\dots,s+k\right\}$, and for any
$\alpha, y > 0$ and $\theta \in \left( 0,\frac{\pi}{2}\right)$ such
that $\sin \theta > c_d y^{-1} \varepsilon$, there exists a measurable
set $\mathcal{B} \subset [0,\infty) \times \mathbb{R}^d \times
\mathbb{S}^{d-1}$, which may depend on $Z_s$, $t$, and
$\left\{ t_j,v_{s+j},\omega_j,i_j\right\}_{j=1}^k$, such that 
\begin{equation}
\begin{aligned}
& \forall T > 0,\\
& \int_0^T \int_{B_{2R}^d} \int_{\mathbb{S}^{d-1}}
\mathbf{1}_{(\tau,v_{s+k+1},\omega_{k+1}) \in \mathcal{B}}
d\omega_{k+1} dv_{s+k+1} d\tau \leq \\
& \qquad \qquad \leq C_{d,s,k} T R^d \left[
\alpha + \frac{y}{\eta T} + C_{d,\alpha} 
\left( \frac{\eta}{R}\right)^{d-1} + C_{d,\alpha} 
\theta^{(d-1)/2} \right]
\end{aligned}
\end{equation}
and
\begin{equation}
\begin{aligned}
& Z_{s,s+k+1} [ Z_s,t+\tau;t_1+\tau,\dots,t_k+\tau,0;
v_{s+1},\dots,v_{s+k},v_{s+k+1};\\
& \qquad \qquad \qquad \qquad \qquad \qquad \qquad
\omega_1,\dots,\omega_k,\omega_{k+1};i_1,\dots,i_k,i_{k+1} ] \\
& \qquad \qquad  \qquad \qquad \qquad \qquad \qquad
\qquad \qquad \qquad \in
\mathcal{G}_{(s+k+1)|m} \cap \hat{\mathcal{U}}_{s+k+1}^\eta
\end{aligned}
\end{equation}
whenever $(\tau,v_{s+k+1},\omega_{k+1}) \in
\left( [0,\infty) \times \mathbb{R}^d \times \mathbb{S}^{d-1}\right)
\backslash \mathcal{B}$.
\end{proposition}
\begin{remark}
The smallness of the error estimate comes from setting
$\eta = \varepsilon^\kappa$ (recall $\kappa \in (0,1)$ is fixed),
$y = \varepsilon^{(1+\kappa)/2}$ and
$\theta \sim \varepsilon^{(1-\kappa)/4}$ to satisfy
the constraint $\sin \theta \geq c_d y^{-1} \varepsilon$.
We regard $\alpha, R$ as fixed as $\varepsilon \rightarrow 0$; it is
then found that $\alpha \rightarrow 0$ and $R\rightarrow \infty$
result in a vanishingly small overall error.
Obviously these choices are not uniquely determined but they
suffice for obtaining the convergence. The implicit dependence
on $\alpha$ could be quantified by writing a quantitative version
of Lemma \ref{lemma:l1}, accounting for the size of the Jacobian
for the mapping $\omega \mapsto u_\omega$ (this is a $(d-1)\times
(d-1)$ determinant)
on the set $\left\{ \omega \cdot v \geq |v| \sin \alpha \right\}$.
\end{remark}
\begin{proof}
Claim \emph{(i)} follows immediately from the definitions of
$\mathcal{G}_{(s+k)|m}$ and $\hat{\mathcal{U}}_{s+k}^\eta$,
 so we turn to Claim
\emph{(ii)}. We begin by deleting creation times where two particles
may be too close to each other; this is dangerous because if particles
are concentrated near the created particle at the time of particle creation
 then we will not be able to
prove that the recollision probability is small. We define
\begin{equation}
\mathcal{B}_I = \left\{
\begin{aligned} 
& (\tau,v_{s+k+1},\omega_{k+1}) \\
& \in [0,\infty) \times \mathbb{R}^d \times \mathbb{S}^{d-1}
\end{aligned}
 \left|
\begin{aligned}
& \exists (x^0,v^0),(x^1,v^1) \in \mathcal{J}_{s+k}
Z_{s+k}^\prime \; \textnormal{such that} \\
& \qquad  (x^0,v^0) \neq (x^1,v^1) \\
& \qquad \textnormal{and }
|(x^0-x^1)-(v^0-v^1)\tau| \leq y
\end{aligned}
\right.
\right\}
\end{equation} 
We can estimate the measure of $\mathcal{B}_I$ because
 there is a known bound on the number of
collisions of $s+k$ hard spheres which is independent of the
initial configuration. \cite{BFK1998} Indeed each pair of distinct points in
$\mathcal{J}_{s+k} Z_{s+k}^\prime$ contributes 
$\mathcal{O}(\eta^{-1} y)$ to the measure of $\mathcal{B}_I$ because
two particles can only stay within a distance $y$ for a time
of order $\eta^{-1} y$; here we are using the time integrals explicitly,
and we are also using the fact that
$Z_{s+k}^\prime \in \hat{\mathcal{U}}_{s+k}^\eta$.
 (Note that we have to account for possible
deletions of particles when estimating the cardinality of
 $\mathcal{J}_{s+k}
Z_{s+k}^\prime$, but there are only $s+k$ deletions and the
dynamics between deletions is the usual hard sphere dynamics;
hence, the total number of collisions is still finite and
quantitatively bounded.) We obtain
\begin{equation}
\begin{aligned}
& \int_0^T \int_{B_{2R}^d} \int_{\mathbb{S}^{d-1}}
\mathbf{1}_{(\tau,v_{s+k+1},\omega_{k+1}) \in \mathcal{B}_I}
d\omega dv_{s+k+1} d\tau \leq \\
& \qquad \qquad \qquad \qquad \qquad \qquad \qquad \qquad \qquad
\leq C_{d,s,k} R^d \eta^{-1} y
\end{aligned}
\end{equation}  
\begin{remark}
$C_{d,s,k}$ may grow rapidly with $s,k$ in accordance with
known bounds on the number of collisions of hard spheres. \cite{BFK1998}
\end{remark}

We define $Z_{s+k}^\prime (\tau) =
\psi_{s+k}^{-\tau} Z_{s+k}^\prime$.
Let us delete particle creations which are too close to grazing:
\begin{equation}
\mathcal{B}_{II} = \left\{
\begin{aligned}
& (\tau,v_{s+k+1},\omega_{k+1}) \in [0,\infty) \times
\mathbb{R}^d \times \mathbb{S}^{d-1} \textnormal{ such that } \\
& \left| \omega_{k+1} \cdot \left(
v_{s+k+1} - v_{i_{k+1}}^\prime (\tau) \right) \right| \leq
(\sin \alpha) \left| v_{s+k+1} - v_{i_{k+1}}^\prime (\tau)\right|
\end{aligned}
\right\}
\end{equation}
We have
\begin{equation}
\int_0^T \int_{B_{2R}^d} \int_{\mathbb{S}^{d-1}}
\mathbf{1}_{(\tau,v_{s+k+1},\omega_{k+1}) \in \mathcal{B}_{II}}
d\omega_{k+1} dv_{s+k+1} d\tau \leq C_d T R^d \alpha
\end{equation}
The remainder of the proof is split between pre-collisional and
post-collisional cases. Here pre-collisional means that
$\omega_{k+1} \cdot \left( v_{s+k+1}-v_{i_{k+1}}^\prime (\tau)\right)
\leq 0$, and
post-collisional means that
$\omega_{k+1} \cdot \left( v_{s+k+1}-v_{i_{k+1}}^\prime (\tau)\right) > 0$.
For convenience we define
\begin{equation}
\mathcal{A}^+ = \left\{ 
\begin{aligned}
& (\tau,v_{s+k+1},\omega_{k+1}) \subset [0,\infty) \times \mathbb{R}^d
\times \mathbb{S}^{d-1} \textnormal{ such that } \\
& \qquad \qquad \qquad \qquad
\omega \cdot \left( v_{s+k+1} - v_{i_{k+1}}^\prime (\tau) \right)
> 0
\end{aligned}
\right\}
\end{equation}
\begin{equation}
\mathcal{A}^- = \left\{ 
\begin{aligned}
& (\tau,v_{s+k+1},\omega_{k+1}) \subset [0,\infty) \times \mathbb{R}^d
\times \mathbb{S}^{d-1} \textnormal{ such that } \\
& \qquad \qquad \qquad \qquad
\omega \cdot \left( v_{s+k+1} - v_{i_{k+1}}^\prime (\tau) \right)
\leq 0
\end{aligned}
\right\}
\end{equation}

\emph{Pre-collisional case.} This is the easier case. We first make
sure that the $(s+k+1)$-particle state is in
$\hat{\mathcal{U}}_{s+k+1}^\eta$ at the time the particle is created.
Note that the existing
particles are at least $y > c_d \varepsilon$ apart at the
time of particle creation.
We define
\begin{equation}
\mathcal{B}_{III}^- =
\left\{
\begin{aligned}
& (\tau,v_{s+k+1},\omega_{k+1}) \in \mathcal{A}^- \backslash
\left( \mathcal{B}_I \cup \mathcal{B}_{II}\right) 
\textnormal{ such that } \\
& \exists (x^0,v^0) \in 
\mathcal{J}_{s+k} \left( Z_{s+k}^\prime (\tau)\right)
\textnormal{ : } 
|v^0 - v_{s+k+1}| \leq \eta
\end{aligned}
\right\} 
\end{equation}
We have
\begin{equation}
\int_0^T \int_{B_{2R}^d} \int_{\mathbb{S}^{d-1}}
\mathbf{1}_{(\tau,v_{s+k+1},\omega_{k+1})\in \mathcal{B}_{III}^-}
d\omega dv_{s+k+1} d\tau \leq C_{d,s,k} T \eta^d
\end{equation}
The constant $C_{d,s,k}$ depends on bounds on the number of
collisions of hard spheres. \cite{BFK1998}

The final estimate is to control recollisions under the backwards
flow. We define
\begin{equation}
\mathcal{B}_{IV}^- =
\left\{
\begin{aligned}
& (\tau,v_{s+k+1},\omega_{k+1}) \in 
\mathcal{A}^- \backslash \left( \mathcal{B}_I \cup \mathcal{B}_{II}\right)
\textnormal{ such that } \\
& \exists (x^0,v^0) \in    \mathcal{J}_{s+k}
\left( Z_{s+k}^\prime (\tau)\right)  \; : \;\\
& \qquad \qquad \frac{\left(x_{i_{k+1}}^\prime (\tau) + \varepsilon \omega
 - x^0 \right)
\cdot \left( v_{s+k+1}-v^0 \right)}
{\left| x_{i_{k+1}}^\prime (\tau) + \varepsilon \omega - x^0\right|
\left|v_{s+k+1}-v^0\right|}
\geq \cos \theta
\end{aligned}
\right\}
\end{equation}
Then we have
\begin{equation}
\int_0^T \int_{B_{2R}^d} \int_{\mathbb{S}^{d-1}}
\mathbf{1}_{(\tau,v_{s+k+1},\omega_{k+1})\in \mathcal{B}_{IV}^-}
d\omega dv_{s+k+1} d\tau \leq C_{d,s,k} T R^d \theta^{d-1}
\end{equation}
As usual the constant $C_{d,s,k}$ could be large.

Let $\mathcal{B}^- = \mathcal{B}_I \cup \mathcal{B}_{II}
\cup \mathcal{B}_{III}^- \cup \mathcal{B}_{IV}^-$; then we have
\begin{equation}
\begin{aligned}
& \int_0^T \int_{B_{2R}^d} \int_{\mathbb{S}^{d-1}}
\mathbf{1}_{(\tau,v_{s+k+1},\omega_{k+1})\in \mathcal{B}^-}
d\omega dv_{s+k+1} d\tau \leq \\
& \qquad \qquad \qquad \qquad \qquad \qquad
\leq C_{d,s,k} T R^d \left[ \alpha + \frac{y}{\eta T} + 
\left(\frac{\eta}{R}\right)^d
+ \theta^{d-1} \right]
\end{aligned}
\end{equation}
But by assumption we have $\sin \theta > c_d y^{-1} \varepsilon$; we
can choose $c_d$ large enough that
\begin{equation}
\begin{aligned}
& Z_{s,s+k+1} [Z_s,t+\tau;t_1+\tau,\dots,t_k+\tau,0;v_{s+1},\dots,
v_{s+k},v_{s+k+1};\\
& \qquad \qquad \qquad \qquad \qquad \qquad \qquad \qquad
\omega_1,\dots,\omega_k,\omega_{k+1},i_1,\dots,i_k,i_{k+1}] \\
& \qquad \qquad \qquad \qquad \qquad \qquad \qquad \qquad \qquad
 \in \mathcal{G}_{(s+k+1)|m} \cap \hat{\mathcal{U}}_{s+k+1}^\eta
\end{aligned}
\end{equation}
whenever $(\tau,v_{s+k+1},\omega_{k+1}) \in \mathcal{A}^- \backslash
\mathcal{B}^-$.

\emph{Post-collisional case.}
Let us define
\begin{equation}
v_{s+k+1}^* = v_{s+k+1} - \omega_{k+1} \omega_{k+1} \cdot
\left( v_{s+k+1} - v_{i_{k+1}}^\prime (\tau)\right)
\end{equation}
\begin{equation}
v_{i_{k+1}}^{\prime *} (\tau) =
v_{i_{k+1}}^\prime (\tau) + \omega_{k+1} \omega_{k+1} \cdot
\left( v_{s+k+1} - v_{i_{k+1}}^\prime (\tau)\right)
\end{equation}
We have to be sure that \emph{both} the $s+k+1$ particle and
the $i_{k+1}$ particle do not recollide with other particles under
the backwards flow. First let us define
\begin{equation}
\mathcal{B}_{III}^+ =
\left\{
\begin{aligned}
& (\tau,v_{s+k+1},\omega_{k+1}) \in \mathcal{A}^+ \backslash
\left( \mathcal{B}_I \cup \mathcal{B}_{II}\right) 
\textnormal{ such that } \\
& \exists (x^0,v^0) \in 
\mathcal{J}_{s+k} \left( Z_{s+k}^\prime (\tau)\right)
\textnormal{ : } 
|v^0 - v_{s+k+1}^*| \leq \eta
\end{aligned}
\right\} 
\end{equation}
\begin{equation}
\mathcal{B}_{IV}^+ =
\left\{
\begin{aligned}
& (\tau,v_{s+k+1},\omega_{k+1}) \in \mathcal{A}^+ \backslash
\left( \mathcal{B}_I \cup \mathcal{B}_{II}\right) 
\textnormal{ such that } \\
& \exists (x^0,v^0) \in 
\mathcal{J}_{s+k} \left( Z_{s+k}^\prime (\tau)\right)
\textnormal{ : } 
|v^0 - v_{i_{k+1}}^{\prime *}(\tau)| \leq \eta
\end{aligned}
\right\} 
\end{equation}
\begin{equation}
\mathcal{B}_{V}^+ =
\left\{
\begin{aligned}
& (\tau,v_{s+k+1},\omega_{k+1}) \in \mathcal{A}^+ \backslash
\left( \mathcal{B}_I \cup \mathcal{B}_{II}\right) 
\textnormal{ such that } \\
& \qquad \qquad \qquad \qquad
|v_{s+k+1} - v_{i_{k+1}}^{\prime}(\tau)| \leq \eta
\end{aligned}
\right\} 
\end{equation}
Then using Lemma \ref{lemma:l1} and Lemma \ref{lemma:l2}, and
bounds on the number of collisions of hard spheres \cite{BFK1998}, we have
\begin{equation}
\int_0^T \int_{B_{2R}^d} \int_{\mathbb{S}^{d-1}}
\mathbf{1}_{(\tau,v_{s+k+1},\omega_{k+1})\in \mathcal{B}_{III}^+}
d\omega dv_{s+k+1} d\tau \leq C_{d,s,k} C_{d,\alpha}
 T R \eta^{d-1}
\end{equation}
\begin{equation}
\int_0^T \int_{B_{2R}^d} \int_{\mathbb{S}^{d-1}}
\mathbf{1}_{(\tau,v_{s+k+1},\omega_{k+1})\in \mathcal{B}_{IV}^+}
d\omega dv_{s+k+1} d\tau \leq C_{d,s,k} C_{d,\alpha}
 T R \eta^{d-1}
\end{equation}
\begin{equation}
\int_0^T \int_{B_{2R}^d} \int_{\mathbb{S}^{d-1}}
\mathbf{1}_{(\tau,v_{s+k+1},\omega_{k+1})\in \mathcal{B}_{V}^+}
d\omega dv_{s+k+1} d\tau \leq C_d
 T \eta^d
\end{equation}

Now we define
\begin{equation}
\mathcal{B}_{VI}^+ =
\left\{
\begin{aligned}
& (\tau,v_{s+k+1},\omega_{k+1}) \in 
\mathcal{A}^+ \backslash \left( \mathcal{B}_I \cup \mathcal{B}_{II}\right)
\textnormal{ such that } \\
& \exists (x^0,v^0) \in    \mathcal{J}_{s+k}
\left( Z_{s+k}^\prime (\tau)\right)  \; : \;
(x^0,v^0) \neq (x_{i_{k+1}}^\prime (\tau),
v_{i_{k+1}}^\prime (\tau)) \\
& \qquad \qquad \textnormal{ and }
 \frac{\left(x_{i_{k+1}}^\prime (\tau) + \varepsilon \omega
 - x^0 \right)
\cdot \left( v_{s+k+1}^*-v^0 \right)}
{\left| x_{i_{k+1}}^\prime (\tau) + \varepsilon \omega - x^0\right|
\left|v_{s+k+1}^*-v^0\right|}
\geq \cos \theta
\end{aligned}
\right\}
\end{equation}
\begin{equation}
\mathcal{B}_{VII}^+ =
\left\{
\begin{aligned}
& (\tau,v_{s+k+1},\omega_{k+1}) \in 
\mathcal{A}^+ \backslash \left( \mathcal{B}_I \cup \mathcal{B}_{II}\right)
\textnormal{ such that } \\
& \exists (x^0,v^0) \in    \mathcal{J}_{s+k}
\left( Z_{s+k}^\prime (\tau)\right)  \; : \;
(x^0,v^0) \neq (x_{i_{k+1}}^\prime (\tau),
v_{i_{k+1}}^\prime (\tau)) \\
& \qquad \qquad \textnormal{ and }
 \frac{\left(x_{i_{k+1}}^\prime (\tau)
 - x^0 \right)
\cdot \left( v_{i_{k+1}}^{\prime *}(\tau)-v^0 \right)}
{\left| x_{i_{k+1}}^\prime (\tau) - x^0\right|
\left|v_{i_{k+1}}^{\prime *}(\tau)-v^0\right|}
\geq \cos \theta
\end{aligned}
\right\}
\end{equation}
Using Lemma \ref{lemma:l1} and Lemma \ref{lemma:l2}, and bounds
on the number of collisions of hard spheres \cite{BFK1998}, we have
\begin{equation}
\int_0^T \int_{B_{2R}^d} \int_{\mathbb{S}^{d-1}}
\mathbf{1}_{(\tau,v_{s+k+1},\omega_{k+1})\in \mathcal{B}_{VI}^+}
d\omega dv_{s+k+1} d\tau \leq C_{d,s,k} C_{d,\alpha} T R^d
 \theta^{(d-1)/2}
\end{equation}
\begin{equation}
\int_0^T \int_{B_{2R}^d} \int_{\mathbb{S}^{d-1}}
\mathbf{1}_{(\tau,v_{s+k+1},\omega_{k+1})\in \mathcal{B}_{VII}^+}
d\omega dv_{s+k+1} d\tau \leq C_{d,s,k} C_{d,\alpha} T R^d
 \theta^{(d-1)/2}
\end{equation}

Let $\mathcal{B}^+ = 
\mathcal{B}_I \cup \mathcal{B}_{II} \cup
\mathcal{B}_{III}^+ \cup \mathcal{B}_{IV}^+ \cup \mathcal{B}_V^+
\cup \mathcal{B}_{VI}^+ \cup \mathcal{B}_{VII}^+$; then we have
\begin{equation}
\begin{aligned}
& \int_0^T \int_{B_{2R}^d} \int_{\mathbb{S}^{d-1}}
\mathbf{1}_{(\tau,v_{s+k+1},\omega_{k+1})\in \mathcal{B}^+}
d\omega dv_{s+k+1} d\tau \leq \\
& \qquad \qquad \qquad \qquad 
\leq C_{d,s,k} T R^d \left[ \alpha + \frac{y}{\eta T} + 
C_{d,\alpha} \left(\frac{\eta}{R}\right)^{d-1}
+ C_{d,\alpha} \theta^{(d-1)/2} \right]
\end{aligned}
\end{equation}
By assumption, $\sin \theta > c_d y^{-1} \varepsilon$; as long
as $c_d$ is chosen large enough, we have
\begin{equation}
\begin{aligned}
& Z_{s,s+k+1} [ Z_s,t+\tau; t_1+\tau,\dots,t_k+\tau,0;
v_{s+1},\dots,v_{s+k},v_{s+k+1};\\
& \qquad \qquad \qquad \qquad \qquad \qquad
\omega_1,\dots,\omega_k,\omega_{k+1};i_1,\dots,i_k,i_{k+1}] \\
& \qquad \qquad \qquad \qquad \qquad \qquad \qquad \qquad
 \in \mathcal{G}_{(s+k+1)|m}\cap \hat{\mathcal{U}}^\eta_{s+k+1}
\end{aligned}
\end{equation}
whenever $(\tau,v_{s+k+1},\omega_{k+1}) \in
\mathcal{A}^+ \backslash \mathcal{B}^+$.
\end{proof}

\section*{Acknowledgements}
\label{sec:acknowledgements}

R.D. gratefully acknowledges support from a postdoctoral fellowship
at the University of Texas at Austin. Additionally, R.D. thanks the
anonymous referee(s) for helpful comments which have 
led to significant improvements in the manuscript.

\bibliography{correlations}

\end{document}